\documentclass[final]{article}
\usepackage{xspace}
\usepackage{amsfonts,amsmath,amssymb,amsthm}
\usepackage{titlesec}
\usepackage{comment}
\usepackage{fancyhdr}
\usepackage[a4paper,body={145mm,210mm},centering]{geometry}
\usepackage[bookmarks,a4paper,colorlinks]{hyperref}
\hypersetup{bookmarksnumbered,plainpages=false,hypertexnames=false}
\hypersetup{pagecolor=blue,linkcolor=blue,citecolor=blue,urlcolor=blue}
\hypersetup{pdfcreator=PDFLaTeX}
\newcommand{\rset}{\mathbf{R}}

\newcommand{\nl}{\nolimits}
\newcommand{\ep}{\varepsilon}
\newcommand{\ind}{\mathbf{1}}
\newcommand{\fl}{\longrightarrow}
\newcommand{\e}{\mathbb{E}}
\newcommand{\p}{\mathbb{P}}
\newcommand{\lp}{\mathrm{L}}
\newcommand{\m}{\mathcal}
\DeclareMathOperator{\sgn}{sgn}
\newcommand{\D}{\mathrm{(D)}}
\newcommand{\ti}[1]{\widetilde{#1}}
\newcommand{\ys}[1][p]{\mathcal{S}^{#1}}
\newcommand{\zs}[1][p]{\mathrm{M}^{#1}}
\newcommand{\dotafter}[1]{#1.}
\titleformat{\section}[hang]
{\normalfont\large\bfseries}{\thesection.}{.5em}{\dotafter}[]
\titleformat{\subsection}[runin]
{\normalfont\bfseries}{\thesubsection.}{.4em}{}[.]
\titlespacing*{\subsection}{0pt}{3ex plus 1ex minus .2ex}{.9em}
\theoremstyle{plain}
\newtheorem{thm}{Theorem}
\newtheorem{lemme}[thm]{Lemma}
\newtheorem{prop}[thm]{Proposition}
\newtheorem{cor}[thm]{Corollary}
\theoremstyle{definition}
\newtheorem*{df}{Definition}
\newtheorem{hyp}{Assumption}

\theoremstyle{remark}
\newtheorem*{rem}{Remark}
\title{\bf Quadratic BSDEs with convex generators and unbounded terminal conditions}
\author{Philippe Briand \and Ying Hu}
\date{\emph{IRMAR, Universit{\'e} Rennes 1, 35 042 Rennes Cedex, FRANCE} \\[2ex]
{\normalsize\tt philippe.briand@univ-rennes1.fr,\quad ying.hu@univ-rennes1.fr}}
\begin{document}
%
%
\maketitle
%
%
\begin{quote}
\textbf{Abstract.}\quad In \cite{BH06}, the authors proved an existence result for BSDEs with quadratic generators with respect to the variable $z$ and with unbounded terminal conditions. However, no uniqueness result was stated in that work. The main goal of this paper is to fill this gap. In order to obtain a comparison theorem for this kind of BSDEs, we assume that the generator is convex with respect to the variable $z$. Under this assumption of convexity, we are also able to prove a stability result in the spirit of the a priori estimates stated in \cite{EPQ97}. With these tools in hands, we can derive the nonlinear Feynman--Kac formula in this context.
\end{quote}
\thispagestyle{fancy}

\section{Introduction}\label{sec:introduction} 

Since the first existence and uniqueness result for backward stochastic differential equations (BSDEs in short for the remaining of the paper) of Pardoux and Peng~\cite{PP90}, lots of works have been done in this area and the original Lipschitz assumption on the generator, the function $f$ in the BSDE
\begin{equation}\label{eq:main}
	Y_t = \xi + \int_t^T f(s,Y_s,Z_s)\, ds - \int_t^T Z_s\cdot dB_s, \qquad 0\leq t\leq T,
\end{equation} 
has been weakened in many situations. Let us recall that, in the previous equation, we are looking for a pair of processes $(Y,Z)$ which is required to be adapted with respect to the filtration generated by the Brownian motion $B$. Even though there exist some improvements of Pardoux--Peng's result  for multidimensional BSDEs (see e.g. \cite{Par99} or \cite{BDHPS03}), the case of real valued BSDEs meaning that $Y$ is a real process is easier to handle. The main reason for the gap between the real case and the multidimensional one comes from the comparison theorem which turns out to be a powerful tool to construct solutions for real valued BSDEs. Roughly speaking, a comparison theorem for BSDEs allows one to compare solutions to BSDEs as soon as one can compare the terminal conditions ($\xi$ in the previous equation) and the generators. In particular, a classical strategy in order to obtain a solution to the BSDE~\eqref{eq:main} consists in constructing an increasing sequence of solutions to BSDEs and then passing to the limit. This method is used for instance in \cite{LsM97} for continuous generators with a linear growth. But one of the most important application of this approach is the work of \textsc{M. Kobylanski}~\cite{Kob00} concerning quadratic BSDEs. We should point out that quadratic BSDE means a BSDE whose generator has a quadratic growth with respect to the variable $z$. For these quadratic BSDEs, all the classical results, existence and uniqueness, comparison theorem and stability of solutions, has been stated in \cite{Kob00} but with the restriction that the terminal conditions have to be bounded random variables.

It was known for a long time that the boundedness of the terminal condition is not really needed for studying a quadratic BSDE. Actually, on simple examples, on can see that the existence of an exponential moment of sufficiently large order is enough to get a solution. In \cite{BH06}, existence of solutions to quadratic BSDEs with unbounded terminal conditions was studied and it was proved that the existence of an exponential moment of $\xi$ is a sufficient condition in order to construct a solution. However the question of uniqueness of solutions in this setting was left open in that work. Quadratic BSDEs with bounded terminal conditions has an interesting feature : $Z \bullet B$, the stochastic integral of the process $Z$ with respect to the Brownian motion $B$, is a BMO--martingale. This observation can be used to obtain uniqueness see \cite{HIM05} and also other properties \cite{BC07P}. When the terminal condition is not bounded, $Z\bullet B$ is, in general, no more a BMO--martingale: this is a difficulty for uniqueness. 

The main contribution of this paper is to obtain a uniqueness result for quadratic BSDEs with unbounded terminal conditions when the generator of \eqref{eq:main} is assumed to be in addition convex or concave with respect to the variable $z$. Let us mention that convex generators appear naturally for BSDEs associated to stochastic control problems,  see e.g. \cite{FHT06}. With this further assumption, it is possible to obtain a comparison theorem for solutions to unbounded quadratic BSDEs which of course gives uniqueness see Section~\ref{sec:comp}. The key idea for proving this result is the following: instead of trying to estimate the difference of two solutions, say $Y$ and $Y'$, we estimate, for each $\theta\in(0,1)$, $Y-\theta Y'$; this allows to take advantage of the convexity of the generator. Moreover, it turns out that the convexity of the generator with respect to the variable $z$ is also a convenient assumption to obtain a stability result. The proof of this last result relies mainly on the same computation mentioned before even though technical difficulties arise and impose us to go into details carefully. Finally all these results yield the nonlinear Feynman--Kac formula in this framework. 

The paper is organized as follows. In Section~\ref{sec:exist}, we prove an existence result in the spirit of \cite{BH06}: here, we work with generators which are quadratic with respect to $z$ and monotone with respect to $y$. This monotonicity assumption allows as usual to get rid of the growth assumption on the generator with respect to $y$. Section~\ref{sec:comp} is devoted to the comparison theorem from which we get as a byproduct an existence and uniqueness result for quadratic BSDEs with unbounded terminal conditions. In Section~\ref{sec:stability}, we obtain a stability result in our framework and finally in the last section we derive the nonlinear Feynman--Kac formula.

\medskip

Let us close this introduction by giving the notations which we will use in all the paper. For the remaining of the paper, let us fix a nonnegative real number $T>0$. First of all, $B=\{B_t\}_{t\geq 0}$ is a standard Brownian
motion with values in $\rset^d$ defined on some complete probability space $(\Omega,\mathcal{F},\p)$. $\{\mathcal{F}_t\}_{t\geq 0}$ is the natural filtration of the Brownian motion $B$ augmented by the $\p$--null sets of $\m F$.  All the measurability notions will refer to this filtration. In particular, the sigma-field of predictable subsets of $[0,T]\times\Omega$ is denoted $\m P$. 

As mentioned in the introduction, we will deal only with real valued BSDEs which is an equation of type~\eqref{eq:main}. The function $f$ is called the generator and $\xi$ the terminal condition.
Let us recall that a generator is  a random function $f:[0,T]\times\Omega\times\rset\times\rset^{d}\fl\rset$ which is measurable with respect to $\m P\otimes\mathcal{B}(\rset)\otimes\mathcal{B}(\rset^d)$ and a terminal condition is simply a real $\mathcal{F}_T$--measurable random variable.

By a solution to the BSDE~\eqref{eq:main} we mean a pair $(Y,Z)=\{ (Y_t,Z_t)\}_{t\in[0,T]}$ of predictable processes  with values in $\rset\times\rset^{d}$ such that $\p$--a.s., $t\longmapsto Y_t$ is continuous, $t\longmapsto Z_t$ belongs to  $\lp^2(0,T)$, $t\longmapsto f(t,Y_t,Z_t)$ belongs to $\lp^1(0,T)$
and $\p$--a.s. 
\begin{equation*}
	Y_t = \xi +\int_t^T f(s,Y_s,Z_s)\, ds -\int_t^T Z_s\cdot dB_s, \qquad 0\leq t\leq T.
\end{equation*} 
We will sometimes use the notation BSDE($\xi,f$) to say that we consider the BSDE whose generator is $f$ and whose terminal condition is $\xi$; $\left(Y^f(\xi),Z^f(\xi) \right)$ means a solution to the BSDE($\xi,f$).
A solution $\left(Y^f(\xi),Z^f(\xi) \right)$ is said to be minimal if $\p$-a.s., for each $t\in[0,T]$, $ Y^f_t(\xi) \leq Y^g_t(\zeta)$ whenever $\p$--a.s. $\xi \leq \zeta$ and $f(t,y,z)\leq g(t,y,z)$ for all $(t,y,z)$. $\left(Y^f(\xi),Z^f(\xi) \right)$ is said to be minimal in some space $\m B$ if it belongs to this space and the previous property holds true as soon as $\left(Y^g(\zeta),Z^g(\zeta) \right)\in\m B$.

For any real $p\geq 1$, $\mathcal{S}^p$ denotes the set of real-valued, adapted and c\`adl\`ag processes
$\{Y_t\}_{t\in[0,T]}$ such that
\begin{equation*}
	\left\| Y \right\|_{\m S^p} :=\e\left[\sup\nl_{0\leq t\leq T} |Y_t|^p\right]^{1/p} < +\infty.
\end{equation*} 
$\left(\ys,\|\cdot\|_{\mathcal{S}^p}\right)$ is a Banach space.

$\mathrm{M}^p\left(\rset^d\right)$ or simply $\zs$ denotes the set of (equivalent classes of) predictable processes $\{Z_t\}_{t\in[0,T]}$ with values in $\rset^d$ such that
\begin{equation*}
	\left\| Z \right\|_{\mathrm{M}^p} : = \e\left[\left(\int_0^T |Z_s|^2\,ds\right)^{p/2}\right]^{1/p} < +\infty.
\end{equation*}
$\mathrm{M}^p\left(\rset^d\right)$ is a Banach space endowed with this norm. We set $\m S = \cup_{p>1} \m S^p$, $\mathrm{M} = \cup_{p>1} \mathrm{M}^p$ and $\m S^\infty$ stands for the set of predictable bounded processes.

Finally, we will say that a real process $Y$, adapted and c\`adl\`ag, belongs to $\m E$ if the random variable $Y^*=\sup_{0\leq t\leq T}|Y_t|$ has exponential moments of all orders and we recall that $Y$ belongs to the class $\D$ as soon as the family $\{ Y_\tau : \tau\leq T\text{ stopping time}\}$ is uniformly integrable.


\section{An existence result}\label{sec:exist} 

In this section, we prove a mere generalization of the existence result for quadratic BSDEs we obtained in \cite{BH06}. We consider here the case where the generator has some monotonicity property with respect to the variable $y$. As usual, this kind of assumptions allows to get rid of the (linear) growth of the generator $f$ in $y$. See e.g. \cite{Par99}, \cite{BDHPS03} or \cite{BLsM06P}. Of course, we still consider the case of a generator with a quadratic growth in the variable $z$ and of an unbounded terminal condition. In this section we assume the following on the generator.

\begin{hyp}\label{H:exist}
	\renewcommand{\labelenumi}{(\roman{enumi})}
	There exist two constants $\beta\geq 0$ and $\gamma >0$ together with a progressively measurable nonnegative stochastic process $\{\alpha(t)\}_{0\leq t\leq T}$ and a deterministic continuous nondecreasing function $\varphi:\rset_+\fl\rset_+$ with $\varphi(0)=0$ such that, $\p$--a.s.,
	\begin{enumerate}
		\item for all $t\in[0,T]$, $(y,z)\longmapsto f(t,y,z)$ is continuous;
		\item monotonicity in $y$ : for each $(t,z)\in[0,T]\times\rset^d$,
		\begin{equation*}
			\forall y\in\rset,\qquad y \left(f(t,y,z)-f(t,0,z)\right) \leq \beta |y|^2;
		\end{equation*}
		\item growth condition:
		\begin{equation*}
			\forall (t,y,z)\in[0,T]\times\rset\times\rset^d,\qquad |f(t,y,z)| \leq \alpha(t) + \varphi(|y|) + \frac{\gamma}{2} |z|^2.
		\end{equation*}
	\end{enumerate}
\end{hyp}

Let us denote by $h$ the (random) function $h(t,x)= \alpha(t) + \beta x$. First of all, we want to derive a sharp estimate for solutions to the BSDE~\eqref{eq:main}. The idea is the following: we apply It{\^o}--Tanaka's formula to compute $U(t,|Y_t|)=e^{\gamma \psi(t,|Y_t|)}$. Let us denote $L$ the local time of $Y$ at 0. We have
\begin{equation*}
	|Y_t| = |\xi| + \int_{t}^T \sgn(Y_s)f(s,Y_s,Z_s)\, ds - \int_t^T \sgn(Y_s) Z_s\cdot dB_s - \int_t^T  dL_s,
\end{equation*}
and, from the growth of $f$ \ref{H:exist}(ii)--(iii),
\begin{align}
	\nonumber
	\sgn(Y_s)f(s,Y_s,Z_s) & = \sgn(Y_s)\left(f(s,Y_s,Z_s)-f(s,0,Z_s)\right) + \sgn(Y_s)f(s,0,Z_s)\\
	\label{eq:majgen}
	&\leq \beta |Y_s| + \alpha(s) + \frac{\gamma}{2} |Z_s|^2	= h(s,|Y_s|) + \frac{\gamma}{2} |Z_s|^2.
\end{align}
From It{\^o}'s formula, we derive the inequality
\begin{align*}
	\frac{dU(t,|Y_t|)}{\gamma U(t,|Y_t|)} & =   \left( -\psi_x(t,|Y_t|) \sgn(Y_t)f(t,Y_t,Z_t) + \psi_t(t,|Y_t|) +\frac{\gamma}{2}\psi_x(t,|Y_t|)^2 |Z_t|^2\right) dt\\
	  & \quad +\frac{1}{2} \psi_{xx}(t,|Y_t|)|Z_t|^2 dt + \psi_x(t,|Y_t|)\, dL_t + \psi_x(t,|Y_t|)\sgn(Y_t)Z_t\cdot dB_t.
\end{align*}
If $\psi_x(t,x)\geq 1$  for $x\geq 0$, we have, taking into account the inequality~\eqref{eq:majgen},
\begin{equation*}
	\psi_x(t,|Y_t|) \sgn(Y_t)f(t,Y_t,Z_t) - \psi_t(t,|Y_t|) -\frac{\gamma}{2}\psi_x(t,|Y_t|)^2 |Z_t|^2 \leq \psi_x(t,|Y_t|) h(t,|Y_t|) - \psi_t(t,|Y_t|).
\end{equation*}
The idea is now clear; we want to find a solution $\psi$ to the PDE: for $s\in[0,T]$, 
\begin{equation*}
	\psi_t(t,x) - h(t,x)\psi_x(t,x) = 0,\qquad \psi(s,x)=\Psi(x),\qquad t\in[s,T], 
\end{equation*}
such that $\psi_x(t,x)\geq 1$  and $\psi_{xx}(t,x)\geq 0$ for $x\geq 0$. But the previous PDE can be solved by the characteristics method. Let $\{v(u;t,x)\}_{0\leq u\leq t}$ be the solution to the integral solution 
\begin{equation*}
	v(u;t,x) = x + \int_u^t h(r,v(r;t,x))\, dr, \quad 0\leq u\leq t.
\end{equation*} 
Then $\psi(t,x) = \Psi\left(v(s;t,x)\right)$. In our case, $h(r,x)= \alpha(r)+\beta x$, the solution is given by the formula
\begin{equation*}
	v(u;t,x) = xe^{\beta(t-u)} + \int_u^t \alpha(r) e^{\beta(r-u)}\, dr.
\end{equation*}
We choose $\Psi(x)=x$ and we get, if $\sigma$ is a stopping time such that $s\leq \sigma \leq T$,
\begin{equation}\label{eq:majimp}
	e^{\gamma |Y_s|} = U(s,|Y_s|)  \leq  U(\sigma,|Y_\sigma|) - \int_s^\sigma \gamma U(t,|Y_t|)\psi_x(t,|Y_t|)\sgn(Y_t)Z_t\cdot dB_t.
\end{equation}

This computation leads the following a priori estimate. 

\begin{prop}\label{en:apest}
	Let $(Y,Z)$ be a solution to the BSDE~\eqref{eq:main} such that the process 
	\begin{equation*}
		\exp\left(\gamma e^{\beta T}|Y_t| + \gamma \int_0^T \alpha(r) e^{\beta r}\,dr\right)
	\end{equation*}  
	belongs to the class $\D$. Then, $\p$--a.s.,
	\begin{equation}\label{eq:estap}
		\forall t\in[0,T],\qquad |Y_t| \leq \frac{1}{\gamma} \log \e\left( \exp\left(\gamma e^{\beta(T-t)} |\xi| + \gamma \int_t^T \alpha(r) e^{\beta (r-t)} dr\right)\:\Big|\: \m F_t \right).
	\end{equation}
\end{prop}

\begin{proof}
	Let $s\in[0,T]$. Set, as before $\psi(t,x)= xe^{\beta(t-s)} + \int_s^t \alpha(r) e^{\beta(r-s)}\, dr$, $U(t,x)=e^{\gamma \psi(t,x)}$. Let us consider, for each integer $n\geq 1$, the following stopping time:
	\begin{equation*}
		\sigma_n = \inf\left\{ t\in[s,T] : \gamma^2 \int_s^t U^2(r,|Y_r|)\psi_x(r,|Y_r|)^2|Z_r|^2\, dr \geq n  \right\}\wedge T.
	\end{equation*}
	It follows from the inequality \eqref{eq:majimp} and the definition of $\sigma_n$ that
	\begin{equation*}
		e^{\gamma |Y_s|} \leq \e\left( e^{\gamma \psi(\sigma_n,|Y_{\sigma_n}|)}\: \big|\: \m F_s\right) \leq \e\left( \exp\left(\gamma e^{\beta(T-s)} |Y_{\sigma_n}| + \gamma \int_s^T \alpha(r)e^{\beta(r-s)}dr\right)\:\Big|\: \m F_s\right).
	\end{equation*}
	Thus, the inequality follows by sending $n$ to infinity.
\end{proof}

\begin{rem}
	It is clear, from the previous computations, that an a priori estimate can be stated for functions $h$ which are not linear as we did. For instance, it is possible to obtain such an estimate when 
	\begin{equation*}
		 y \left(f(t,y,z)-f(t,0,z)\right) \leq |y| \rho(|y|)
	\end{equation*}  
	where $\rho: \rset_+\fl\rset_+$ is convex, $\m C^1$ with $\rho(0)=0$; $|f(t,y,z)| \leq \alpha + \varphi(|y|) + \frac{\gamma}{2} |z|^2$ ($\alpha>0$) and  
	\begin{equation*}
		\int_0^{+\infty} \frac{du}{h(u)} = +\infty,\quad\text{with }h(u)=\rho(u)+\alpha.
	\end{equation*}
	In this case, one has to choose $\psi(t,x)= \Theta^{-1}(t-s + \Theta(x))$ where 
	\begin{equation*}
		\Theta(x) = \int_0^x \frac{du}{h(u)}.
	\end{equation*}
	See e.g. \cite{LsM98} or \cite{BH06}.
\end{rem}

As a byproduct of this a priori estimate, we see that, when $|\xi|$ and $|\alpha|_1:=\int_0^T \alpha(s)\, ds$ are bounded random variables, any solution $(Y,Z)$ to \eqref{eq:main} such that $e^{\gamma e^{\beta T} |Y_t|}$ belongs to the class $\D$, is actually bounded. More precisely,
\begin{equation*}
	|Y_t| \leq e^{\beta T} \, \| |\xi| + |\alpha|_1 \|_\infty.
\end{equation*}
This observation allows to generalize a little bit, with exactly the same proofs, some of Kobylanski's results \cite{Kob00}. 

\begin{lemme}[M. Kobylanski]\label{en:kob}
	Let the assumption \ref{H:exist} hold. If $\xi$ and $|\alpha|_1$ are bounded random variables, then the BSDE~\eqref{eq:main} has a minimal bounded solution (and a maximal one also).

	Moreover, let $(f_n)_{n\geq 1}$ be a sequence of generators satisfying \ref{H:exist} with $(\alpha,\beta,\gamma,\varphi)$ independent of $n$ such that $\left(f_n(t,y_n,z_n)\right)_{n\geq 1}$ converges to $f(t,y,z)$ as soon as $(y_n,z_n)\fl (y,z)$ and let $(\xi_n)_{n\geq 1}$ be a sequence of terminal conditions converging almost surely to $\xi$. Let us assume that, for each $n\geq 1$, the BSDE associated to $\xi_n$ and $f_n$ has a solution $(Y^n,Z^n)\in\ys[\infty]\times\zs[2]$ such that $(Y^n)_{n\geq 1}$ is nondecreasing and, for some constant $C$, $\sup_{n\geq 1}\| Y^n\|_\infty \leq C$. 
	
	Then, there exists $(Y,Z)\in\ys[\infty]\times\zs[2]$, solution to the BSDE~\eqref{eq:main}, such that $(Y^n)_{n\geq 1}$ converges to $Y$ uniformly on $[0,T]$ in probability and $(Z^n)_{n\geq 1}$ converges to $Z$ in $\zs[2]$.
\end{lemme}

We are now in position to generalize to this framework the result we obtained in \cite{BH06}.

\begin{prop}\label{en:exist}
	Let \ref{H:exist} hold.
	
	If $|\xi|+|\alpha|_1$ has an exponential moment of order $\gamma e^{\beta T}$, then the BSDE~\eqref{eq:main} has a solution $(Y,Z)$ such that 
	\begin{equation}\label{eq:majY}
		\forall t\in[0,T],\qquad |Y_t| \leq \frac{1}{\gamma} \log \e\left( \exp\left(\gamma e^{\beta(T-t)} |\xi| + \gamma \int_t^T \alpha(r) e^{\beta (r-t)} dr\right)\:\Big|\: \m F_t \right).
	\end{equation}	
\end{prop}

\begin{proof}
	Let us assume that $\xi$ and $f$ are nonnegative. For each integer $n\geq 1$, we consider $\xi_n=\xi\wedge n $ and $f_n(t,y,z) = \ind_{t\leq \sigma_n} f(t,y,z)$ where
	\begin{equation*}
		\sigma_n = \inf\left\{t\in[0,T] : \int_0^t \alpha(s)\, ds \geq n \right\} \wedge T.
	\end{equation*}
According to the first part of the Lemma~\ref{en:kob}, let $(Y^n,Z^n)\in\ys[\infty]\times\zs[2]$ be the minimal bounded solution to the BSDE whose terminal condition is $\xi_n$ and whose generator is $f_n$. It follows, from Proposition~\ref{en:apest}, that for each $n\geq 1$,
\begin{equation*}
	\left| Y^n_t \right| \leq \frac{1}{\gamma}\, \log \e\left( \exp\left(\gamma e^{\beta T} \left(|\xi| + |\alpha|_1\right)\right)\:|\: \m F_t\right) := X_t.
\end{equation*} 
Since $\xi_n\leq \xi_{n+1}$, $f_n\leq f_{n+1}$ and since we are dealing with minimal solutions, the sequence $(Y^n)$ is nondecreasing. We have in mind to use the stability property quoted in the second part of Lemma~\ref{en:kob}. To use this result, we need to know that the sequence $(Y^n)$ remains bounded which, of course, is not true in general. To overcome this difficulty the idea is to use the previous estimate and to work on random time interval instead of working on the whole interval $[0,T]$. To do this, let, for $k\geq 1$, 
$\tau_k$ be the following stopping time
\begin{equation*}
	\tau_k = \inf\{t\in[0,T] : X_t\geq k\}\wedge T.
\end{equation*}
By construction the sequence (in $n$) $Y^n_k(t)=Y^n_{t\wedge \tau_k}$ remains bounded in $k$. Setting moreover $Z^n_k(t)=Z^n_t\, \ind_{t\leq \tau_k}$, we have
\begin{equation*}
	Y^n_k(t) = Y^n_{\tau_k} + \int_t^T \ind_{t\leq \tau_k\wedge\sigma_n} f\left(s,Y^n_k(s),Z^n_k(s)\right) ds - \int_t^T Z^n_k(s)\cdot dB_s.
\end{equation*}
We apply the second part of Lemma~\ref{en:kob} to obtain, for each $k$, a solution $(Y_k,Z_k)$ to the BSDE
\begin{equation*}
	Y_k(t) = \xi_k + \int_t^{\tau_k} f\left(s,Y_k(s),Z_k(s)\right) ds - \int_t^{\tau_k} Z_k(s)\cdot dB_s,\quad\text{ with } \xi_k = \sup\nl_{n\geq 1} Y^n_{\tau_k}.
\end{equation*}
But by construction, $\tau_k\leq \tau_{k+1}$ so we have the localization property
\begin{equation*}
	Y_{k+1}(t\wedge\tau_k) = Y_k(t),\qquad Z_{k+1}(t)\ind_{t\leq \tau_k} = Z_k(t).
\end{equation*}
If we define the processes $Y$ and $Z$ on $[0,T]$ by the formula
\begin{equation*}
	Y_t = Y_1(t)\ind_{t\leq \tau_1} + \sum_{k\geq 2} Y_k(t)\,\ind_{]\tau_{k-1},\tau_k]}(t),\qquad Z_t = Z_1(t)\ind_{t\leq \tau_1} + \sum_{k\geq 2} Z_k(t)\,\ind_{]\tau_{k-1},\tau_k]}(t),
\end{equation*}
the previous BSDE is rewritten as follows
\begin{equation*}
	Y_t = \xi_k + \int_t^{\tau_k}f\left(s,Y_s,Z_s\right) ds - \int_t^{\tau_k} Z_s\cdot dB_s.
\end{equation*}
The last point is the fact that, $\p$--a.s., $\tau_k=T$ for $k$ large enough which allows to send $k$ to $+\infty$ in the previous equation to prove that $(Y,Z)$ is a solution to \eqref{eq:main}. Of course, the inequality of the theorem is verified by the process $Y$ since it is verified by each process $Y^n$ in view of Proposition~\ref{en:apest}.

In the general case, we use a double approximation; $\xi^{n,p}=\xi^+\wedge n - \xi^-\wedge p$,
\begin{equation*}
	f^{n,p}(t,y,z) = \ind_{t\leq \sigma_n}f^+(t,y,z) - \ind_{t\leq\sigma_p} f^-(t,y,z).
\end{equation*}
\end{proof}

\begin{cor}\label{en:corint}
	Let the assumption \ref{H:exist} hold. 
	
	If $|\xi| + |\alpha|_1$ has an exponential moment of order $\lambda>\gamma e^{\beta T}$ then the BSDE~\eqref{eq:main} has a solution $(Y,Z)$ such that $Y$ satisfies the inequality~\eqref{eq:majY} and $Z$ belongs to $\zs[2]$. 
	
	 In the case where $|\xi| + |\alpha|_1$ has exponential moments of all order, this solution is such that $Y$ belongs to $\m E$ and $Z$ belongs to $\zs$ for all $p\geq 1$. More precisely, for each $p>1$, 
	\begin{equation}\label{eq:MajAp}
		\e\left[\exp\left(\gamma p \sup_{0\leq t \leq T}|Y_t|\right) + \left(\int_0^T |Z_s|^2 \, ds\right)^{p/2}\right]  \leq C \, \e\left[ \exp\bigl(p\gamma \left(|\xi| + |\alpha|_1\right)\bigr) \right]
	\end{equation}
	where $C$ depends on $p$, $\gamma$, $\beta$ and $T$.
\end{cor}

\begin{proof}
	Let $(Y,Z)$ the solution to \eqref{eq:main} obtained in Theorem~\ref{en:exist}. Then $Y$ satisfies the estimate \eqref{eq:majY}.  It follows immediately from this estimate and Doob's maximal inequality that $e^{\gamma Y^*}$ belongs to $\ys[q]$ for some $q>1$ when $\zeta:=|\xi| + |\alpha|_1$ has an exponential moment of order $\lambda>\gamma e^{\beta T}$ and that the estimate~\eqref{eq:MajAp} holds true for $Y$ when $\zeta$ has exponential moments of all order. To obtain the results for the process $Z$, we argue as in \cite{BH06}. 
	
	For $n\geq 1$, let $\tau_n$ be the following stopping time
	$$
	\tau_n = \inf\left\{ t\geq 0 : \int_0^t e^{2\gamma |Y_s|} |Z_s|^2 \, ds \geq n \right\} \wedge T,
	$$
	and let us consider the function from $\rset_+$ into itself defined by
	$$
	u(x) = \frac{1}{\gamma^2} \left( e^{\gamma x} - 1 - \gamma x \right).
	$$
	$x\longmapsto u(|x|)$ is $\m C^2$ and we have from It{\^o}'s formula, with the notation $\sgn(x) = - \ind_{x\leq 0} + \ind_{x>0}$,
	\begin{align*}
	u(|Y_0|) & =  u(|Y_{t\wedge\tau_n}|) + \int_0^{t\wedge\tau_n} \left( u'(|Y_s|) \sgn(Y_s) f(s,Y_s,Z_s) - \frac{1}{2} u''(|Y_s|) |Z_s|^2 \right) ds  \\
	&  \quad - \int_0^{t\wedge\tau_n} u'(|Y_s|) \sgn(Y_s) Z_s\cdot dB_s.
	\end{align*}
	It follows from \ref{H:exist} that
	\begin{equation*}
		\sgn(Y_s) f(s,Y_s,Z_s) \leq \sgn(Y_s) \left(f(s,Y_s,Z_s)-f(s, 0,Z_s)\right) + |f(s,0,Z_s)| \leq \alpha(s) + \beta |Y_s| + \frac{\gamma}{2} |Z_s|^2,
	\end{equation*} 
	and since $u'(x) \geq 0$ for $x\geq 0$ that
	\begin{align*}
	u(|Y_0|) & \leq  u(|Y_{t\wedge\tau_n}|) + \int_0^{t\wedge\tau_n} u'(|Y_s|) \left(\alpha(s)  + \beta |Y_s|\right) ds - \int_0^{t\wedge\tau_n} u'(|Y_s|) \sgn(Y_s) Z_s\cdot dB_s \\
	& \quad - \frac{1}{2} \int_0^{t\wedge\tau_n} \left(u''(|Y_s|) - \gamma\, u'(|Y_s|) \right) |Z_s|^2 ds .
	\end{align*}
	Moreover, we have $(u'' - \gamma u')(x) = 1$ and $u(x)\geq 0$ for $x\geq 0$, so taking the expectation of the previous inequality
	\begin{equation*}
		\frac{1}{2}\, \e\left[\int_0^{\tau_n} |Z_s|^2\,ds\right] \leq \e\left[\frac{1}{\gamma^2}\,\sup_{t\in[0,T]} e^{\gamma |Y_t|} + \frac{1}{\gamma}\, \int_0^T  e^{\gamma |Y_s|} \left( \alpha(s)  + \beta |Y_s| \right) ds\right].
	\end{equation*} 
	Fatou's lemma together with the fact that $e^{\gamma |Y_t|}\in \m S^q$ for some $q>1$ gives the result when $\zeta$ has an exponential moment of order $\lambda > \gamma e^{\beta T}$.
	
	For the second part of the result, let us observe that
	\begin{align*}
		\frac{1}{2}\,\int_0^{\tau_n} |Z_s|^2\,ds & \leq \frac{1}{\gamma^2}\, \sup_{0\leq t\leq T} e^{\gamma |Y_t|} + \frac{1}{\gamma}\, \int_0^T  e^{\gamma |Y_s|} \left( \alpha(s)  + \beta |Y_s| \right) ds  \\
		& \quad  + \sup_{0\leq t\leq T}\left\vert \int_0^{t\wedge\tau_n} u'(|Y_s|) \sgn(Y_s) Z_s\cdot dB_s\right\vert.
	\end{align*}
	It follows from the BDG inequalities, since $|y| \leq e^{\gamma|y|}/\gamma$ and $u'(|y|)\leq e^{\gamma |y|}/\gamma$, that
	\begin{align*}
		\e\left[\left(\int_0^{\tau_n} |Z_s|^2\,ds\right)^{p/2}\right] & \leq C\,\e\left[e^{p\gamma Y^*} + e^{p\gamma Y^*/2} |\alpha|_1^{p/2} + \left(\int_0^{\tau_n} e^{2\gamma|Y_s| }|Z_s|^2\, ds\right)^{p/4}\right] \\
		& \leq C'\, \e\left[e^{p\gamma Y^*} + e^{p\gamma Y^*/2} |\alpha|_1^{p/2}\right] + \frac{1}{2}\,	\e\left[\left(\int_0^{\tau_n} |Z_s|^2\,ds\right)^{p/2}\right]
	\end{align*}
from which the result follows using Fatou's lemma.	
\end{proof}


\section{Comparison Theorem}\label{sec:comp} 

Let us consider now the main topic of this paper : uniqueness for quadratic BSDEs with unbounded terminal condition. We state in this section a comparison theorem for solutions to quadratic BSDEs and let us recall first that, in the bounded case, such a result exists (see \cite{Kob00}) even though it requires more assumptions than those in the existence result. In our unbounded framework the situation is exactly the same : the assumptions to prove uniqueness are stronger than those needed for existence. In particular, a very convenient way to derive the comparison theorem is to assume that the generator is convex (or concave) with respect to the variable $z$. From the point of view of PDEs, the convexity of the generator seems a natural assumption see e.g. \cite{dLL06}. Let us consider the following assumption on the generator $f$.

\begin{hyp}\label{H:comp}
	\renewcommand{\labelenumi}{(\roman{enumi})}
	There exist two constants $\gamma >0$ and $\beta\geq 0$ together with a nonnegative progressively measurable stochastic process $\{\alpha(t)\}_{0\leq t\leq T}$ such that, $\p$--a.s.,
	\begin{enumerate}
		\item for all $t\in[0,T]$, for all $y\in\rset$, $z\longmapsto f(t,y,z)$ is convex;
		\item for all $(t,z)\in[0,T]\times\rset^d $, 
		\begin{equation*}
			\forall \left(y,y'\right)\in\rset^2,\qquad \left|f(t,y,z)-f(t,y',z)\right| \leq \beta \left|y-y'\right| ;
		\end{equation*}
		\item $f$ has the following growth:
		\begin{equation*}
			\forall (t,y,z)\in[0,T]\times\rset\times\rset^d,\qquad |f(t,y,z)| \leq \alpha(t) + \beta\,|y| + \frac{\gamma}{2}\, |z|^2 ;
		\end{equation*}
		\item $|\alpha|_1$ has exponential moment of all order.
	\end{enumerate}
\end{hyp} 

\pagebreak[2]

\begin{thm}[Comparison theorem]\label{en:comp}
Let $(Y,Z)$ be a solution to \eqref{eq:main} and $(Y',Z')$ be a solution to the BSDE associated to the terminal condition $\xi'$ and to the generator $f'$ such that both $Y$ and $Y'$ belongs to $\m E$, namely, 
\begin{equation*}
	\forall \lambda >0, \qquad \e\left[ e^{\lambda Y^*} + e^{\lambda {Y'}^*}\right] <+\infty.
\end{equation*}
We assume that, $\p$--a.s.,
\begin{equation*}
	\xi\leq \xi',\qquad \forall(t,y,z)\in[0,T]\times\rset\times\rset^d,\quad f(t,y,z) \leq f'(t,y,z).
\end{equation*} 

If $f$ verifies \ref{H:comp}, then $\p$--a.s., for each $t\in[0,T]$, $Y_t\leq Y'_t$.

If moreover, $Y_0=Y'_0$, then 
\begin{equation*}
	\p\left(\xi'-\xi = 0,\ \int_0^T \left(f'-f\right)\left(t,Y'_t,Z'_t\right) dt = 0\right) >0.
\end{equation*}
\end{thm}

\begin{proof}
The idea is the following: instead of trying to estimate the difference between the processes $Y$ and $Y'$, we estimate $Y-\theta Y'$, for each $\theta\in(0,1)$, in order to take advantage of the convexity of the generator. Similar idea is also used for quadratic PDEs in \cite{dLL06} where the authors prove that $u-\theta u'$ is nonpositive for two solutions $u$ and $u'$.
	
Let $\theta\in(0,1)$ and let us set  $U_t = Y_t- \theta Y'_t$ and $V_t = Z_t -\theta Z'_t$. Let us consider a real stochastic process $\{a(t)\}_{0\leq t\leq T}$, progressively measurable, with integrable paths to be chosen later. We set, for all $t\in[0,T]$, $A_t=\int_0^t a(s)\, ds$. We have, from It{\^o}'s formula,
$$
e^{A_t}U_t = e^{A_T}U_T + \int_t^T e^{A_s}F_s\, ds -\int_t^T e^{A_s} V_s\cdot dB_s,\quad 0\leq t\leq T,
$$
where, of course,
\begin{equation}\label{eq:dfF}
	F_t=  \left(f(t,Y_t,Z_t)-\theta f'\left(t,Y'_t,Z'_t\right)\right) - a(t) U_t.
\end{equation}
Let us set as usual  $\delta f(t) = \left(f-f'\right)\left(t,Y'_t,Z'_t\right)$ so that  
\begin{multline*}
	  f(t,Y_t,Z_t)-\theta f'\left(t,Y'_t,Z'_t\right)  = f(t,Y_t,Z_t)-\theta f\left(t,Y'_t,Z'_t\right) + \theta \delta f(t)   \\
	 =\left(f(t,Y_t,Z_t)- f\left(t,Y'_t,Z_t\right)\right) + \left(f\left(t,Y'_t,Z_t\right) - \theta f\left(t,Y'_t,Z'_t\right)\right) + \theta \delta f(t).
\end{multline*}
 
Since $f$ is convex with respect to $z$, the second term of the right hand side of the previous inequality can be easily bounded from above. Indeed,
\begin{align*} 
	f\left(t,Y'_t,Z_t\right) & =  f\left(t,Y'_t, \theta Z'_t + (1-\theta)\frac{Z_t-\theta Z'_t}{1-\theta}\right) \\
 	&\leq   \theta f\left(t,Y'_t,Z'_t\right) + (1-\theta)  f\left(t,Y'_t,\frac{Z_t-\theta Z'_t}{1-\theta}\right) 
\end{align*}
and from the growth of the generator $f$, \ref{H:comp}(iii),
\begin{equation}\label{eq:maj1}
	f\left(t,Y'_t,Z_t\right)\leq  \theta f\left(t,Y'_t,Z'_t\right) + (1-\theta) \left(\alpha(t)+\beta \left|Y'_t\right|\right) +\frac{\gamma}{2(1-\theta)}\left|Z_t-\theta Z'_t\right|^2.
\end{equation}

Roughly speaking, the first term does not change anything since it depends more or less only on $y$ and $f$ is Lipschitz with respect to this variable. We get rid of it with a classical linearization. Let us write
\begin{align*}
	f\left(t,Y_t,Z_t\right)-f\left(t,Y'_t,Z_t\right) & = f\left(t,Y_t,Z_t\right)-f\left(t,\theta Y'_t,Z_t\right) + f\left(t,\theta Y'_t,Z_t\right)-f\left(t,Y'_t,Z_t\right) \\
	& = a(t) U_t + f\left(t,\theta Y'_t,Z_t\right)-f\left(t,Y'_t,Z_t\right) 
\end{align*}
where $a(t) = \left[f\left(t,Y_t,Z_t\right)-f\left(t,\theta Y'_t,Z_t\right)\right] / U_t$ when $U_t\neq 0$ and $a(t)=\beta$ in the other case. Since $f$ is $\beta$--Lipschitz, $a$ is bounded by $\beta$ and
\begin{equation}\label{eq:maj2}
	f\left(t,Y_t,Z_t\right)-f\left(t,Y'_t,Z_t\right) \leq a(t) U_t + (1-\theta)\beta \left|Y'_t\right|.
\end{equation} 

If we choose for $a$ the process we have just introduced, we get, from \eqref{eq:maj1} and \eqref{eq:maj2}, coming back to the definition of $F$ \eqref{eq:dfF},
\begin{equation}\label{eq:majF}
	F_t \leq  (1-\theta) \left(\alpha(t)+ 2\beta \left|Y'_t\right|\right)+\frac{\gamma}{2(1-\theta)}|V_t|^2  +\theta\delta f(t)  .
\end{equation}

Now, we get rid of the quadratic term with an exponential change of variables. Let $c\geq 0$ and let us set $P_t = e^{c e^{A_t}U_t}$, $Q_t = c  P_t  V_t e^{A_t}$. From It{\^o}'s formula we deduce that
\begin{align*}
P_t & =  P_T + c\int_t^T P_se^{A_s}  \left(  F_s  - \frac{ce^{A_s}}{2} \left|V_s\right|^2 \right) ds  - c\int_t^T e^{A_s}P_sV_s\cdot dB_s \\
& := P_T + \int_t^T G_s\, ds - \int_t^T Q_s\cdot dB_s.
\end{align*}
\eqref{eq:majF} yields, since $c$ is nonnegative,
\begin{equation*}
	G_t \leq cP_te^{A_t} \left\{ (1-\theta) \left( \alpha(t)+ 2\beta \left|Y'_t\right| \right) +\theta\delta f(t)\right\} +cP_te^{A_t} \left(\frac{\gamma }{2(1-\theta)} -\frac{ce^{A_t}}{2}\right) \left| V_t\right|^2.
\end{equation*}
But $A_t \geq -\beta T$ so that, if we choose $c=c(\theta):=\gamma e^{\beta T}/(1-\theta)$ we obtain the following inequality, 
\begin{equation}\label{eq:forstab}
	G_t \leq P_te^{A_t} \left( \theta c(\theta)\delta f(t) + \gamma e^{\beta T} \left(\alpha(t)+ 2 \beta \left|Y'_t\right|\right) \right).
\end{equation}

Finally, let us introduce the processes
\begin{equation*}
	D_t = \exp \left(\int_0^t e^{A_s}\left( \theta c(\theta) \delta f(s) + \gamma e^{\beta T}  \left(\alpha(s)+2\beta \left|Y'_s\right|\right) \right) ds\right), \quad \ti P_t = D_t P_t,\quad \ti Q_t = D_t Q_t.
\end{equation*}
Once again It{\^o}'s formula gives us, for any stopping time $\tau$ such that $0\leq t\leq \tau \leq T$,
\begin{equation*}
	 \ti P_t \leq \ti P_\tau - \int_t^\tau \ti Q_s\cdot dB_s.
\end{equation*}
Let us consider, for $n\geq 1$, $\tau_n$ the stopping time
\begin{equation*}
	\tau_n = \inf\left\{ u \geq t : \int_t^u \left| \ti Q_s \right|^2 ds \geq n \right\} \wedge T.
\end{equation*}
We get from the previous equation
\begin{equation*}
	P_t \leq \e\left( \exp \left(\int_t^{\tau_n} e^{A_s}\left( \theta c(\theta)\delta f(s) + \gamma e^{\beta T}  \left(\alpha(s)+2\beta \left|Y'_s\right|\right) \right) ds\right)  P_{\tau_n} \:\Big|\: \m F_t\right), 
\end{equation*}
and, in view of the integrability assumption on $\alpha$, $Y$ and $Y'$, since $|A_s|\leq \beta T$, we can send $n$ to infinity to obtain
\begin{equation*}
	 P_t \leq \e\left( \exp \left(\int_t^T e^{A_s} \left( \theta c(\theta)\delta f(s) + \gamma e^{\beta T}  \left(\alpha(s)+2\beta \left|Y'_s\right|\right) \right) ds\right)  P_{T} \:\Big|\: \m F_t\right) .
\end{equation*}
On the other hand, $\left(\xi-\theta\xi'\right) = (1-\theta) \xi + \theta \left(\xi-\xi'\right) \leq (1-\theta) |\xi| + \theta\delta\xi$, and since $\delta\xi$ and $\delta f$ are nonpositive, we finally derive the inequality
\begin{multline}\label{eq:majSC}
 	  \exp\left(\frac{\gamma e^{\beta T + A_t}}{1-\theta}\left(Y_t-\theta Y'_t\right)\right)  \\ 
	 \leq  \e\left( \exp\left\{ \frac{\gamma\theta}{1-\theta}\left( \delta\xi + \int_t^T   \delta f(s)\, ds\right) + \gamma e^{2\beta T} \left(|\xi| + \int_t^T \left(\alpha(s)+2\beta \left|Y'_s\right|\right)  ds\right)  \right\} \:\bigg|\: \m F_t\right). 
\end{multline}    

In particular, since $\delta\xi$ and $\delta f$ are nonpositive and $\beta T + A_t \geq 0$,
\begin{equation*}
	Y_t-\theta Y'_t \leq \frac{1-\theta}{\gamma}\, \log \e\left( \exp\left\{ \gamma e^{2\beta T} \left(|\xi| + \int_t^T \left(\alpha(s)+2\beta \left|Y'_s\right|\right)  ds\right)  \right\} \:\bigg|\: \m F_t\right),
\end{equation*}
and sending $\theta$ to 1, we get $Y_t-Y'_t\leq 0$ which gives the first part of the result.

For the second part of the theorem, if in addition we have $Y_0=Y'_0$, then, coming back to the inequality~\eqref{eq:majSC}, we deduce that
\begin{equation*}
	e^{\gamma e^{\beta T} Y_0} \leq \e\left[ \exp\left\{ \frac{\gamma\theta}{1-\theta}\left( \delta\xi + \int_0^T   \delta f(s)\, ds\right) + \gamma e^{2\beta T} \left(|\xi| + \int_0^T \left(\alpha(s)+2\beta \left|Y'_s\right|\right)  ds\right)  \right\} \right],
\end{equation*}
and, once again, sending $\theta$ to 1, we get
\begin{equation*}
	0 < e^{\gamma e^{\beta T} Y_0} \leq \e\left[ \ind_{ \delta\xi + \int_0^T \delta f(s)\, ds = 0}\exp\left\{ \gamma e^{2\beta T} \left(|\xi| + \int_0^T \left(\alpha(s)+2\beta \left|Y'_s\right|\right)  ds\right)  \right\} \right]
\end{equation*}
which says that $\p\left(\delta \xi =0, \int_0^T \delta f(s)\, ds =0\right) >0$.
\end{proof}

As a byproduct, we obtain, using Corollary~\ref{en:corint}, the following existence and uniqueness result.

\begin{cor}\label{en:exun}
	Let the assumption~\ref{H:comp} hold and let us assume moreover that the random variables $|\xi|$ and $|\alpha|_1$ have exponential moments of all order. Then the BSDE~\eqref{eq:main} has a unique solution $(Y,Z)$ such that $Y$ belongs to $\m E$ and $Z$ belongs to $\zs$ for each $p\geq 1$. 
\end{cor}


\section{Stability}\label{sec:stability} 

As we have seen in the previous section, the convexity of the generator with respect to the variable $z$ leads to the comparison theorem. It turns out that we can also derive a stability result under this assumption. To be more precise, let us consider a generator $f$ for which \ref{H:comp} holds with parameters $(\alpha,\beta,\gamma)$ together with a sequence of generators $(f_n)_{n\geq 1}$ : for each $n\geq 1$, $f_n$ satisfy \ref{H:comp} with parameters $(\alpha_n,\beta,\gamma)$. Finally, let $\xi$ and $(\xi_n)_{n\geq 1}$ be random terminal values such that, for each $\lambda >0$,
\begin{equation}\label{eq:conv}
	\e\left[e^{\lambda\left(|\xi|+|\alpha|_1\right)}\right] + \sup\nl_{n\geq 1}\e\left[e^{\lambda\left(|\xi_n|+|\alpha_n|_1\right)}\right] < +\infty.
\end{equation}

According to Corollary~\ref{en:exun}, let $(Y,Z)$ be the solution to the BSDE~\eqref{eq:main} and let us introduce, for each $n\geq 1$, $(Y^n,Z^n)$ the solution to the BSDE
\begin{equation*}
	Y^n_t = \xi_n + \int_t^T f_n(s,Y^n_s,Z^n_s)\, ds - \int_t^T Z^n_s\cdot dB_s.
\end{equation*}

\begin{prop}\label{en:stab}
	Let \ref{H:comp} hold for $f$ and $f_n$ and let us assume moreover that the inequality~\eqref{eq:conv} holds true. 
	
	If $\xi_n\fl\xi$ $\p$--a.s. and, $m$ being the Lebesgue measure on $[0,T]$, $m\otimes\p$--a.e., for each $(y,z)\in\rset\times\rset^d$, $f_n(t,y,z)\fl f(t,y,z)$, then, for each $p\geq 1$,
	\begin{equation*}
		\e\left[\exp\left(\sup\nl_{0\leq t\leq T} |Y^n_t-Y_t|\right)^p + \left(\int_0^T |Z^n_s-Z_s|^2\, ds\right)^{p/2}\right]\fl 0.
	\end{equation*}
	
\end{prop}     

\begin{proof}
	It follows from Corollary~\ref{en:corint} and the integrability assumptions~\eqref{eq:conv} that the sequence $\left((Y^n,Z^n)\right)_{n\geq 1}$ satisfies, for each $p\geq 1$,
	\begin{equation*}
		\sup_{n\geq 1} \e\left[ \exp\left(\sup\nl_{0\leq t\leq T} |Y^n_t|\right)^p + \left(\int_0^T |Z^n_s|^2\, ds\right)^{p/2}\right] < +\infty.
	\end{equation*}
	It is thus enough to prove that 
	\begin{equation*}
		\sup\nl_{0\leq t\leq T} |Y^n_t-Y_t| + \int_0^T |Z^n_s-Z_s|^2\, ds 
	\end{equation*}
	converges to 0 in probability to get the result of the proposition.
	
	Let us fix $\theta\in(0,1)$ and $n\geq 1$. First of all, since $f_n$ is convex in $z$ and $\beta$--Lipschitz in $y$, we can argue exactly as in the proof of the comparison theorem (Theorem~\ref{en:comp}) to estimate $Y^n_t-\theta Y_t$. Setting $P_t = e^{ce^{A_t} U_t}$ and $Q_t= cP_te^{A_t} (Z^n_t-\theta Z_t)$ where
	\begin{equation*}
		a_t = \left[f_n\left(t,Y^n_t,Z^n_t\right) - f_n\left(t,\theta Y_t, Z^n_t\right)\right]/ \left[Y^n_t-\theta Y_t\right],\quad A_t = \int_0^t a(s)\, ds, \quad c= \frac{\gamma e^{\beta T}}{1-\theta},
	\end{equation*}
	we have $d P_t = -G_t\, dt + Q_t\cdot dB_t$, with, see the inequality~\eqref{eq:forstab},
	\begin{equation*}
		G_t \leq P_te^{A_t} \left( c\theta \delta_n f(t) + \gamma e^{\beta T} \left(\alpha_n(t)+ 2 \beta \left|Y_t\right|\right) \right), 
	\end{equation*}
	where $\delta_n f(t)$ stands for $\left(f_n-f\right)(t,Y_t,Z_t)$. Buy taking into account the fact that $a(t)$ is bounded by $\beta$, we have
	\begin{equation*}
		G_t \leq P_t \gamma e^{2\beta T} \left( \frac{|\delta_n f(t)|}{1-\theta} + \alpha_n(t) + 2\beta |Y_t|\right).
	\end{equation*}
	Let us introduce as before the process 
	\begin{equation*}
		D^n_t = \exp \left(\gamma e^{2\beta T}\, \int_0^t \left(\alpha_n(s) + 2\beta |Y_s|\right) ds\right);
	\end{equation*}
	we get from It{\^o}'s formula
	\begin{equation*}
		D^n_t P_t \leq \e\left( D^n_T P_T + \frac{\gamma e^{2\beta T}}{1-\theta}\int_t^T P_s D^n_s |\delta_n f(s)| ds\:\Big|\: \m F_t \right) .
	\end{equation*}
  
We observe that, since $a$ is bounded by $\beta$,  
\begin{align*}
   P_s & \leq  G_n(\theta) : = \sup_{0\leq t\leq T}\exp\left( \frac{\gamma e^{2\beta T}}{1-\theta} (|Y_t| + |Y^n_t| )\right) , \\
   P_T & \leq X_n(\theta) : = \exp\left( \frac{\gamma e^{2\beta T}}{1-\theta} (|\xi^n-\theta \xi| \vee  |\xi-\theta \xi^n|)\right).
\end{align*} 

With these notations, we derive the inequality
\begin{equation*}
   Y^n_t-\theta\, Y_t \leq  \frac{(1-\theta) e^{-\beta T -A_t}}{\gamma}  \log  \e\left( D^n_T X_n(\theta) + \frac{\gamma e^{2\beta T}}{1-\theta}\, D^n_T G_n(\theta)\, \int_t^T |\delta_n f(s)| ds\:\Big|\: \m F_t \right)
\end{equation*} 
and finally, since $\log x \leq x$,
\begin{equation}\label{eq:cayest}
	Y^n_t-Y_t \leq (1-\theta) |Y_t|  +  \frac{1-\theta}{\gamma} \e\left(D^n_T X_n(\theta)\:|\: \m F_t\right) + e^{2\beta T} \e\left(D^n_T G_n(\theta)\, \int_t^T |\delta_n f(s)| ds\:\Big|\: \m F_t\right).	
\end{equation}

Now we want to find an upper bound for $Y_t-Y^n_t$ and to do this we use the same approach. Once again, let $\theta\in(0,1)$ and let us set $U_t = Y_t -\theta Y^n_t$, $V_t = Z_t-\theta Z^n_t$. We have
$$
dU_t = -H_t\, dt + V_t\cdot dB_t
$$
with $H_t = f(t,Y_t,Z_t)-\theta f_n(t,Y^n_t,Z^n_t)$. To get the same inequality, we split $H$ in the following way:
\begin{align*}
	H_t & =  f(t,Y_t,Z_t)-f_n(t,Y_t,Z_t) + f_n(t,Y_t,Z_t) - \theta f_n(t,Y^n_t,Z^n_t) \\
 		 & =  -\delta_n f(t) + f_n(t,Y_t,Z_t)-\theta f_n(t,Y_t,Z^n_t) + \theta f_n(t,Y_t,Z^n_t)-\theta f_n(t,Y^n_t,Z^n_t).
\end{align*}
As before, the convexity of $f_n$ with respect to $z$ leads to  
\begin{equation*}
   f_n(t,Y_t,Z_t)  -  \theta f_n(t,Y_t,Z^n_t) \leq (1-\theta) (\alpha_n(t) + \beta |Y_t| ) + \frac{\gamma}{2(1-\theta)} |V_t|^2. 
\end{equation*}
For the third term, let us  introduce the process  
\begin{equation*}
   a(t) = \frac{ f_n(t,Y_t,Z^n_t)-f_n(t,Y^n_t,Z^n_t) }{Y_t-Y^n_t}\ind_{|Y_t-Y^n_t|>0} 
\end{equation*}
which is bounded by $\beta$ so that
\begin{align*}
   \theta f_n(t,Y_t,Z^n_t)-\theta f_n(t,Y^n_t,Z^n_t) & =  \theta a(t) (Y_t-Y^n_t) = a(t) (\theta Y_t -Y_t + Y_t -\theta Y^n_t) \\
&\leq  \beta(1-\theta) |Y_t| + a(t) U_t.
\end{align*}
It follows from the previous inequalities that 
\begin{equation*}
	H_s \leq |\delta_n f(s) | + (1-\theta) (\alpha_n(s) + 2\beta|Y_s| ) + \frac{\gamma}{2(1-\theta)} |V_s|^2 + a(s) U_s
\end{equation*} 
with $a$ bounded by $\beta$. It follows from It{\^o}'s formula that, taking as usual $A_t = \int_0^t a(s)\,ds$, 
\begin{equation*}
	e^{A_t}U_t = e^{A_T}U_T + \int_t^T e^{A_s}F_s\, ds -\int_t^T e^{A_s} V_s\, dB_s,\quad 0\leq t\leq T,
\end{equation*}   
with
\begin{equation*}
	F_s \leq |\delta_n f(s) | + (1-\theta) (\alpha_n(s) + 2\beta|Y_s| ) + \frac{\gamma}{2(1-\theta)} |V_s|^2.
\end{equation*} 
This is the same inequality as~\eqref{eq:majF}. As a by product, we deduce that the inequality~\eqref{eq:cayest} is valid also for $Y_t-Y^n_t$, namely,
\begin{equation*}
   Y_t-Y^n_t \leq (1-\theta) |Y^n_t|  +  \frac{1-\theta}{\gamma} \e\left(D^n_T X_n(\theta)\:|\: \m F_t\right) + e^{2\beta T} \e\left(D^n_T G_n(\theta)\, \int_t^T |\delta_n f(s)| ds\:\Big|\: \m F_t\right)
\end{equation*} 
and finally we have
\begin{multline}\label{eq:enfin}  
	|Y^n_t-Y_t| \leq  (1-\theta) \left(|Y^n_t| + |Y_t|\right)  +  \frac{1-\theta}{\gamma} \e\left(D^n_T X_n(\theta)\:|\: \m F_t\right)  \\ 
	+ e^{2\beta T} \e\left(D^n_T G_n(\theta)\, \int_0^T |\delta_n f(s)| ds\:\Big|\: \m F_t\right).   
\end{multline}  

Let us fix $\ep >0$. We have from~\eqref{eq:enfin} and Doob's maximal inequality
\begin{multline*}
	\p\left(\sup_{0\leq t\leq T} |Y^n_t-Y_t| > \ep\right) \leq  \frac{3(1-\theta)}{\ep} \e\left[\sup_{0\leq t\leq T} \left(|Y^n_t| + |Y_t|\right)\right]  \\
	+ \frac{3(1-\theta)}{\gamma\ep} \e\left[D^n_T X_n(\theta)\right] + e^{2\beta T} \frac{3}{\ep} \e\left[D^n_T G_n(\theta)\, \int_0^T |\delta_n f(s)| ds\right]
\end{multline*}  
and since, for $\theta\in(0,1)$, the sequences $\left(\sup_{0\leq t\leq T} \left(|Y^n_t| + |Y_t|\right)_{n\geq 1}\right)$, $\left(D^n_T\right)_{n\geq 1}$ and $\left(G_n(\theta)\right)_{n\geq 1}$ are bounded in all $\lp^p$ spaces, we deduce from H{\"o}lder's inequality 
\begin{equation}\label{eq:ad}
	\p\left(\sup_{0\leq t\leq T} |Y^n_t-Y_t| > \ep\right) \leq \frac{1-\theta}{\ep} C\left(1+ \|X_n(\theta)\|_2\right) + \frac{C(\theta)}{\ep} \left\| \int_0^T |\delta_n f(s)| ds \right\|_2.
\end{equation}
Let us recall that $X_n(\theta) = \exp\left( \frac{\gamma e^{2\beta T}}{1-\theta} (|\xi^n-\theta \xi| \vee  |\xi-\theta \xi^n|) \right)$ so, as $n$ goes to $\infty$, $X_n(\theta)$ converges to $\exp\left( \gamma e^{2\beta T} |\xi| \right)$ almost surely and actually in all $\lp^p$ spaces in view of the integrability assumptions on the sequence $(\xi_n)_{n\geq 1}$ (see~\eqref{eq:conv}). Moreover, $|\delta_n f(s)| = |f-f_n|(s,Y_s,Z_s)$ converges to 0 $m\otimes\p$--a.e. and, since 
\begin{equation*}
	|\delta_n f(s)| \leq \alpha(s) + \alpha_n(s) + 2\beta |Y_s| + \gamma |Z_s|^2,
\end{equation*}
Corollary~\ref{en:exun}  and the inequality~\eqref{eq:conv} ensures the convergence to 0 of $\int_0^T |\delta_n f(s)|\, ds$ in all $\lp^p$ spaces. The inequality~\eqref{eq:ad} yields, for all $0< \theta < 1$,
\begin{equation*}
	\limsup_{n\to +\infty}\p\left(\sup_{0\leq t\leq T} |Y^n_t-Y_t| > \ep\right) \leq \frac{1-\theta}{\ep} C \left(1+ \left\| \exp\left( \gamma e^{2\beta T} |\xi| \right) \right\|_2\right),
\end{equation*}  
and, sending $\theta$ to 1, we obtain the convergence of $\sup_{0\leq t\leq T}|Y^n_t-Y_t|$ to 0 in probability. As we mentioned before, we obtain the convergence of $e^{\sup_{0\leq t\leq T}|Y^n_t-Y_t|}$ in all $\lp^p$ spaces.

To get the convergence of $\int_0^T |Z^n_s-Z_s|^2 \, ds$ to 0 in probability, let us just mention that, from It{\^o}'s formula we have
\begin{equation*}
	\e\left[\int_0^T  |Z^n_s-Z_s|^2 \, ds \right] \leq \e\left[|\xi_n-\xi|^2 + 2 \sup_{0\leq t\leq T}|Y^n_t-Y_t| \, \int_0^T |f_n(s,Y^n_s,Z^n_s)-f(s,Y_s,Z_s)| \, ds\right]
\end{equation*}
from which the result follows directly.
\end{proof}

\section{Application to quadratic PDEs}\label{sec:application_to_quadratic_pdes} 
In this section, we give an application of our results concerning BSDEs to PDEs which are quadratic with respect to the gradient of the solution. More precisely, we want to obtain the nonlinear Feynman--Kac formula in this framework. Let us consider the following semilinear PDE
\begin{equation}\label{eq:pde}
	\partial_t u(t,x) + \m L u(t,x) + f\left(t,x,u(t,x),\sigma^{*}\nabla_x u(t,x)\right)=0,\qquad u(T,\cdot)=g,
\end{equation}
where $\m L$ is the infinitesimal generator of the diffusion $X^{t_0,x_0}$ solution to the SDE
\begin{equation}\label{eq:pdesde}
	X_t = x_0 + \int_{t_0}^t b(s,X_s)\, ds + \int_{t_0}^t \sigma(s,X_s)\, dB_s, \quad t_0\leq t\leq T,\qquad X_t=x_0,\quad t\leq t_0.
\end{equation}
The nonlinear Feynman--Kac consists in proving that the function defined by the formula
\begin{equation}\label{eq:fk}
	\forall(t,x)\in[0,T]\times\rset^n,\qquad u(t,x):= Y^{t,x}_t
\end{equation}
where, for each $(t_0,x_0)\in[0,T]\times\rset^n$, $\left(Y^{t_0,x_0},Z^{t_0,x_0}\right)$ stands for the solution to the following BSDE
\begin{equation}\label{eq:pdebsde}
	Y_t = g\left(X^{t_0,x_0}_T\right) + \int_t^T f\left(s,X^{t_0,x_0}_s,Y_s,Z_s\right) ds -\int_t^T Z_s\cdot dB_s,\qquad 0\leq t\leq T,
\end{equation}
is a solution, at least a viscosity solution, to the PDE~\eqref{eq:pde}. Our objective is to derive this probabilistic representation for the solution to the PDE when the nonlinearity $f$ is quadratic with respect to $\nabla_x u$ and when $g$ is an unbounded function. 

Let us first give our assumptions concerning the linear part of the PDE namely the coefficients of the diffusion.

\begin{hyp}\label{H:edsquad}
	\renewcommand{\labelenumi}{(\roman{enumi})}
	Let $b: [0,T]\times\rset^n\fl \rset^n$ and $\sigma: [0,T]\times\rset^n\fl \rset^{n\times d}$ be continuous functions and let us assume that there exists $\beta \geq 0$ such that:
	\begin{enumerate}
		\item for all $t\in[0,T]$, $|b(t,0)|\leq \beta$, and
		\begin{equation*}
			\forall \left(x,x'\right)\in\rset^n\times\rset^n,\qquad \left|b(t,x)-b\left(t,x'\right)\right| + \left|\sigma(t,x)-\sigma\left(t,x'\right)\right| \leq \beta \left| x-x'\right| ;
		\end{equation*}
		\item $\sigma$ is bounded.
	\end{enumerate}
\end{hyp}

Under the assumption~\ref{H:edsquad}, for each $(t_0,x_0)\in[0,T]\times\rset^n$, the SDE~\eqref{eq:pdesde} has a unique solution denoted $X^{t_0,x_0}$. Classical results on SDEs show that, for each $p\geq 1$, $X^{t_0,x_0}$ belongs to $\ys$.  Actually, since $\sigma$ is assumed to be a bounded function, for $1\leq p <2$, we have 
\begin{equation*}
\forall \lambda >0,\qquad	\e\left[\sup\nl_{0\leq t\leq T} e^{\lambda \left|X^{t_0,x_0}_t\right|^p}\right] \leq C\, e^{\lambda C |x|^p},
\end{equation*}
where the constant $C$ depends upon $p$, $T$, $\beta$, $\lambda$ and $\|\sigma\|_\infty$. Indeed, we have
\begin{equation*}
	\sup_{t_0\leq t\leq u} \left| X^{t_0,x_0}_t \right| \leq |x_0| + \beta T + \beta \int_{t_0}^u \sup_{t_0\leq t\leq s} \left| X^{t_0,x_0}_t \right| \, ds + \sup_{t_0\leq t\leq T} \left| \int_{t_0}^t \sigma\left(s,X^{t_0,x_0}_s\right) dB_s\right|
\end{equation*}
and we deduce from Gronwall's lemma the inequality
\begin{equation*}
	\sup_{t_0\leq t\leq u} \left| X^{t_0,x_0}_t \right| \leq \left(|x_0| + \beta T + \sup_{t_0\leq t\leq T} \left| \int_{t_0}^t \sigma\left(s,X^{t_0,x_0}_s\right) dB_s\right|\right) e^{\beta T}.
\end{equation*}
It follows from the Dambis--Dubins--Schwarz representation of the continuous martingale 
\begin{equation*}
	\int_0^t \sigma(s)\, dB_s, \quad\text{ with}\quad \sigma(s) =  \sigma\left(s,X^{t_0,x_0}_s\right) \ind_{t_0< s\leq T} + \ind_{s > T}
\end{equation*}
that, for $\lambda >0$ and $1\leq p <2$, 
\begin{equation*}
	\e\left[\sup_{t_0\leq t\leq T}\exp\left(\lambda  \left| \int_{t_0}^t \sigma\left(s,X^{t_0,x_0}_s\right) dB_s\right|^p\right)\right] \leq  \e\left[ \sup_{0\leq t\leq \|\sigma\|_\infty^2 T}e^{\lambda  |B_t|^p}\right],
\end{equation*}
which is a finite constant depending on $p$, $\lambda$, $T$ and $\|\sigma\|_\infty$.

With this observation in hands, we can give our assumptions on the nonlinear term of the PDE, the generator $f$, and the terminal condition.

\begin{hyp}\label{H:edsrquad}
	\renewcommand{\labelenumi}{(\roman{enumi})}
	Let $f: [0,T]\times\rset^n\times\rset\times\rset^{d}\fl \rset$ and $g: \rset^n\fl \rset$ be continuous and let us assume moreover that there exist two constants $\beta \geq 0$ and $1\leq p <2$ such that:
	\begin{enumerate}
		\item for each $(t,x,z)\in[0,T]\times\rset^n\times\rset^d$,
		\begin{equation*}
			\forall \left(y,y'\right)\in\rset^2,\qquad \left|f\left(t,x,y,z\right)-f\left(t,x,y',z\right)\right| \leq \beta \left|y-y'\right|~;
		\end{equation*}
		\item for each $(t,x,y)\in[0,T]\times\rset^n\times\rset$, $z\longmapsto f(t,x,y,z)$ is convex on $\rset^d$;
		\item for each $(t,x,y,z)\in[0,T]\times\rset^n\times\rset\times\rset^d$,
		\begin{equation*}
			|f(t,x,y,z)| + |g(x)| \leq \beta \left(1 + |x|^p + |y| + |z|^2\right).
		\end{equation*}
	\end{enumerate}
\end{hyp}

Since, for $1\leq p<2$, $\sup \left|X^{t_0,x_0}_t\right|^p$ has exponential moments of all orders, the growth condition on $f$ and $g$ allows one to use Corollary~\ref{en:exun} to construct a unique solution, $\left(Y^{t_0,x_0},Z^{t_0,x_0}\right)$, to the BSDE~\eqref{eq:pdebsde}. Moreover, $Y^{t_0,x_0}$ belongs to $\m E$ and the process $Z^{t_0,x_0}$ belongs to the space $\zs$ for each $p\geq 1$. 

It is worth noticing that $Y^{t_0,x_0}_{t_0}$ is actually deterministic for each point $(t_0,x_0)\in[0,T]\times\rset^n$. Indeed, the process $\left\{X^{t_0,x_0}_{t}\right\}_{t_0\leq t\leq T}$ is known to be measurable with respect to the filtration generated by the increments of the Brownian motion after time $t_0$ and it is by now well known that the process $\left\{\left(Y^{t_0,x_0}_t,Z^{t_0,x_0}_t\right)\right\}_{t_0\leq t\leq T}$ inherits this property. As a byproduct of this observation, $u$ defined by the formula~\eqref{eq:fk} is a deterministic function. Let us prove that $u$ is a viscosity solution to the PDE~\eqref{eq:pde}.

\begin{prop}\label{en:quadpde}
	Let the assumption \ref{H:edsquad} and \ref{H:edsrquad} hold.
	
	 The function $u$ defined by \eqref{eq:fk} is continuous on $[0,T]\times\rset^n$ and satisfies 
	\begin{equation*}
		\forall (t,x)\in[0,T]\times\rset^n,\qquad |u(t,x)| \leq C\,(1+|x|^p).
	\end{equation*}
	
	Moreover $u$ is a viscosity solution to \eqref{eq:pde}.
\end{prop}

Before proving this result, let us recall what is a viscosity solution to \eqref{eq:pde}.

\begin{df}
    \label{def:visc}
A continuous function $u$ on $[0,T]\times\rset^n$ such that $u(T,x)=g(x)$ is said to be a viscosity subsolution (respectively supersolution) to \eqref{eq:pde} if
\begin{equation*}
	\partial_t\varphi(t_0,x_0) + \mathcal{L}\varphi(t_0,x_0) + f\left(t_0,x_0, u(t_0,x_0), 
	\sigma^*\nabla_x\varphi(t_0,x_0)\right) \geq 0, \qquad \text{ (respectively 
	}\leq 0\text{)}
\end{equation*} 
as soon as  $u-\varphi$ has a local maximum (respectively minimum) at $(t_0,x_0)\in(0,T)\times\rset^n$ where $\varphi$ is a smooth function. 

A viscosity solution is both a viscosity subsolution and a viscosity supersolution.
\end{df}

\begin{proof}
	Our proof uses standard arguments since, in our setting, the BSDEs have all the properties required: comparison and stability. First of all, the continuity of $u$ is an immediate consequence of the stability property (Proposition~\ref{en:stab}) since the map $(t,x)\longmapsto X^{t,x}$ is known to be continuous. Secondly, the growth of the function $u$ comes directly from the general estimate on $Y$ see \eqref{eq:majY}. 
	
	Let us briefly explain why $u$ is a viscosity subsolution to \eqref{eq:pde}. Let $\varphi$ be a smooth function such that $u-\varphi$ has a local maximum at the point $(t_0,x_0)\in(0,T)\times\rset^n$. Without loss of generality, we can assume that $\varphi(t_0,x_0) = u(t_0,x_0)$ and we want to prove that
	\begin{equation*}
		\partial_t \varphi(t_0,x_0) + \mathcal{L}\varphi(t_0,x_0) +f\left(t_0,x_0,u(t_0,x_0),\sigma^*\nabla_x\varphi(t_0,x_0)\right) \geq 0.
	\end{equation*}
	Let us assume that the previous inequality does not hold. By continuity, there exist $\delta >0$ and $0< \alpha \leq T-t_0$ such that
	\begin{equation*}
		u(t,x) \leq \varphi(t,x),\qquad\text{et,}\qquad \partial_t\varphi(t,x) + \mathcal{L}\varphi(t,x) + 
		f\left(t,x,u(t,x), \sigma^*\nabla_x\varphi(t,x)\right) \leq -\delta
	\end{equation*} 
	as soon as $t_0\leq t \leq t_0+\alpha$ and $|x-x_0| \leq \alpha$.
	
	Let $\tau$ be the following stopping time
	\begin{equation*}
		\tau = \inf\left\{ u\geq t_0 : \left|X^{t_0,x_0}_u -x_0\right| \geq \alpha \right\} \wedge (t_0 + \alpha).
	\end{equation*}
	 The proof consists in applying the comparison theorem, actually the strict version, to the processes
	\begin{equation*}
		\left(Y^{t_0,x_0}_{t\wedge\tau}, \ind_{t\leq\tau} Z^{t_0,x_0}_t\right)
		\quad\text{ and }\quad 
		\left(\varphi\left(t\wedge\tau,X^{t_0,x_0}_{t\wedge\tau}\right),\ind_{t\leq \tau}\sigma^*\nabla_x\varphi\left(t,X^{t_0,x_0}_{t}\right)\right)
	\end{equation*} 
	respectively solution to the BSDEs
	\begin{align*}
		Y_t & = Y_{t_0+\alpha} +\int_t^{t_0+\alpha} \ind_{s\leq\tau}\,f\left(s,X^{t_0,x_0}_s,Y_s, Z_s\right)\, ds - \int_t^{t_0+\alpha} Z_s\cdot dB_s, \\
		Y'_t & = \varphi\left(\tau, X^{t_0,x_0}_\tau\right) + \int_t^{t_0+\alpha} -\ind_{s\leq \tau} 
		\left\{\partial_t\varphi+\mathcal{L}\varphi\right\}\left(s,X^{t_0,x_0}_s\right) \, 
		ds - \int_t^{t_0+\alpha} Z'_s\cdot dB_s.
	\end{align*}
	
	In order to compare the terminal conditions and the generators of these BSDEs, let us recall that the uniqueness of solutions to \eqref{eq:pdebsde} yields the Markov property:
	\begin{equation*}
		Y^{t_0,x_0}_{t}= u\left(t, X^{t_0,x_0}_{t}\right),\quad t\geq t_0.
	\end{equation*} 
	Hence, we can rewrite the first BSDE in the following way
	\begin{equation*}
		Y_t =  u\left(\tau,X^{t_0,x_0}_\tau\right) + \int_t^{t_0+\alpha} 
		\ind_{s\leq\tau}\,f\left(s,X^{t_0,x_0}_s,u\left(s,X^{t_0,x_0}_s\right), Z_s\right)\, ds - \int_t^{t_0+\alpha} Z_s\cdot dB_s.
	\end{equation*}
	By definition of the stopping time $\tau$, we have $u\left(\tau, X^{t_0,x_0}_\tau\right) \leq \varphi\left(\tau,X^{t_0,x_0}_\tau\right)$ together with
	\begin{align*}
	\ind_{s\leq\tau}\,f\left(s,X^{t_0,x_0}_s,u\left(s,X^{t_0,x_0}_s\right), 
	Z'_s\right) & =  
	\ind_{s\leq\tau}\,f\left(s,X^{t_0,x_0}_s,u\left(s,X^{t_0,x_0}_s\right), \sigma^*\nabla_x\varphi(s,X^{t_0,x_0}_s)\right) \\
	& \leq  -\ind_{s\leq \tau} 
	\left\{\partial_t\varphi+\mathcal{L}\varphi\right\}\left(s,X^{t_0,x_0}_s\right) - \delta\ind_{s\leq \tau}.
	\end{align*}
	Moreover, it is worth noticing that
	\begin{equation*}
		\int_{t_0}^{t_0+\alpha}-\ind_{s\leq 
		\tau}\left(\partial_t\varphi+\mathcal{L}\varphi+f\right) \left(s, X^{t_0,x_0}_s, 
		u\left(s,X^{t_0,x_0}_s\right), 
		\sigma^*\nabla_x\varphi\left(s,X^{t_0,x_0}_s\right)\right)\,ds \geq \delta (\tau-t_0)
	\end{equation*}
	and of course $\p(\tau=t_0)=0$. Thus, we can apply the second part of the comparison theorem, Theorem~\ref{en:comp} above, and we get $u(t_0,x_0)= Y_{t_0} < Y'_{t_0}=\varphi(t_0,x_0)$. But we have assumed that $u(t_0,x_0)=\varphi(t_0,x_0)$: $u$  has to be a subsolution to~\eqref{eq:pde}.  The fact that $u$ is a supersolution and thus a solution can be shown in the same way.	
\end{proof}

\begin{rem}
	When $f$ does not depend on the variable $y$ it can be shown that $u$ is the unique viscosity solution with quadratic growth: $|u(t,x)|\leq C\left(1+|x|^2\right)$. This follows from the uniqueness results in \cite{dLL06} concerning Bellman--Isaacs equation. 
\end{rem}


%
%
\providecommand{\bysame}{\leavevmode\hbox to3em{\hrulefill}\thinspace}
\providecommand{\MR}{\relax\ifhmode\unskip\space\fi MR }
\providecommand{\MRhref}[2]{%
  \href{http://www.ams.org/mathscinet-getitem?mr=#1}{#2}
}
\providecommand{\href}[2]{#2}


\begin{thebibliography}{10}

\bibitem{BC07P}
Ph. Briand and F.~Confortola, \emph{{BSDE}s with stochastic {L}ipschitz
  condition and quadratic {PDE}s in {H}ilbert spaces}, Preprint, {\small\url{http://fr.arxiv.org/abs/math.PR/0701849}}, 2007.

\bibitem{BDHPS03}
Ph. Briand, B.~Delyon, Y.~Hu, {\'E}.~Pardoux, and L.~Stoica, \emph{L$^p$
  solutions of backward stochastic differential equations}, Stochastic Process.
  Appl. \textbf{108} (2003), no.~1, 109--129.

\bibitem{BH06}
Ph. Briand and Y.~Hu, \emph{{BSDE} with quadratic growth and unbounded terminal
  value}, Probab. Theory Related Fields \textbf{136} (2006), no.~4, 604--618.

\bibitem{BLsM06P}
Ph. Briand, J.-P. Lepeltier, and J.~San~Mart{\'\i}n, \emph{One-dimensional
  {B}{S}{D}{E}'s whose coefficient is monotonic in y and non-lipschitz in z},
  Bernoulli, to appear, {\small\url{http://perso.univ-rennes1.fr/philippe.briand/docs/blsm.pdf}}, 2006.

\bibitem{dLL06}
F.~Da~Lio and O.~Ley, \emph{Uniqueness results for second-order
  {B}ellman-{I}saacs equations under quadratic growth assumptions and
  applications}, SIAM J. Control Optim. \textbf{45} (2006), no.~1, 74--106
  (electronic).

\bibitem{EPQ97}
N.~El~Karoui, S.~Peng, and M.-C. Quenez, \emph{Backward stochastic differential
  equations in finance}, Math. Finance \textbf{7} (1997), no.~1, 1--71.

\bibitem{FHT06}
M.~Fuhrman, Y.~Hu, and G.~Tessitore, \emph{On a class of stochastic optimal
  control problems related to {BSDE}s with quadratic growth}, SIAM J. Control
  Optim. \textbf{45} (2006), no.~4, 1279--1296 (electronic).

\bibitem{HIM05}
Y.~Hu, P.~Imkeller, and M.~M{\"u}ller, \emph{Utility maximization in incomplete
  markets}, Ann. Appl. Probab. \textbf{15} (2005), no.~3, 1691--1712.

\bibitem{Kob00}
M.~Kobylanski, \emph{Backward stochastic differential equations and partial
  differential equations with quadratic growth}, Ann. Probab. \textbf{28}
  (2000), no.~2, 558--602.

\bibitem{LsM97}
J.-P. Lepeltier and J.~San~Mart{\'\i}n, \emph{Backward stochastic differential
  equations with continuous coefficients}, Statist. Probab. Lett. \textbf{32}
  (1997), no.~4, 425--430.

\bibitem{LsM98}
\bysame, \emph{Existence for {B}{S}{D}{E} with superlinear-quadratic
  coefficient}, Stochastics Stochastics Rep. \textbf{63} (1998), no.~3-4,
  227--240.

\bibitem{Par99}
{\'E}.~Pardoux, \emph{B{S}{D}{E}s, weak convergence and homogenization of
  semilinear {P}{D}{E}s}, Nonlinear analysis, differential equations and
  control (Montreal, QC, 1998), Kluwer Acad. Publ., Dordrecht, 1999,
  pp.~503--549.

\bibitem{PP90}
{\'E}.~Pardoux and S.~Peng, \emph{Adapted solution of a backward stochastic
  differential equation}, Systems Control Lett. \textbf{14} (1990), no.~1,
  55--61.

\end{thebibliography}
\end{document}